\documentclass[preprint,11pt]{elsarticle}

\usepackage{amsmath, amsthm, amssymb,mathtools}
\usepackage{tikz}
\providecommand{\LP}{{\mathcal{P}}}
\providecommand{\LB}{\mathcal{B}}
\newtheorem{thm}{Theorem}
\newtheorem{lem}[thm]{Lemma}
\newtheorem{cor}[thm]{Corollary}
\newtheorem{pro}[thm]{Proposition}
\newtheorem{cons}[thm]{Construction}
\newdefinition{rmk}{Remark}
\newproof{pf}{Proof}
\newproof{pot3}{Proof of Theorem \ref{upbound}}
\newproof{pot4}{Proof of Theorem \ref{asym}}
\newproof{pot5}{Proof of Theorem \ref{upbound2}}
\begin{document}

\begin{frontmatter}



\title{Suitable sets of permutations, packings of triples, and Ramsey's theorem}


\author{Xiande Zhang}

\address{School of Mathematical Sciences,
 University of Science and Technology of China,\\
Wu Wen-Tsun Key Laboratory of Mathematics, No. 96 Jinzhai Road, Hefei, 230026, Anhui, China ({\tt drzhangx@ustc.edu.cn}).}

\begin{abstract}
A set of $N$ permutations of $\{1,2,\ldots,v\}$ is $t$-suitable, if each symbol precedes each subset of $t-1$ others in at least one permutation. The extremal problem of determining the smallest size $N$ of such sets for given $v$ and $t$ was the subject of classical studies by Dushnik in 1950 and Spencer in 1971. Colbourn recently introduced the concept of suitable cores as equivalent objects of suitable sets of permutations, and studied the dual problem of determining the largest $v=\text{SCN}(t,N)$ such that a suitable core exists for given $t$ and $N$. Chan and Jedwab showed that when $N=\lfloor\frac{t+1}{2}\rfloor\lceil\frac{t+1}{2}\rceil+l$, the value of SCN$(t,N)$ is asymptotically $\lfloor\frac{t}{2}\rfloor+2$ if $l$ is a fixed integer. In this paper, we improve  this result by showing that it is also true when $l=O(\ln t)$ using Ramsey theory. When $v$ is bigger than $\lfloor\frac{t}{2}\rfloor+2$,  we give new explicit constructions of suitable cores  from packings of triples, and random constructions from extended Ramsey colorings.

\end{abstract}

\begin{keyword}
Extremal problems \sep Ramsey's theorem \sep Suitable arrays \sep Suitable cores


\end{keyword}

\end{frontmatter}


\section{Introduction}
\label{se:intro}

A set $\LP$ of permutations $\{\pi_1,\ldots,\pi_N\}$ on $[v]=\{1,2,\ldots,v\}$ is called {\it suitable of strength $t$}, or $t$-suitable, if for every subset $S\subset [v]$ of size $t$ and every $\sigma\in S$, there is a permutation $\pi\in \LP$ for which $\pi^{-1}(\sigma)<\pi^{-1}(s)$ for every $s\in S\setminus\{\sigma\}$. Forming an $N\times v$ array $A$ in which the entry in position $(i,j)$  is $\pi_i(j)$, one can
equivalently say that  each  symbol of $[v]$
precedes  each  subset  of  $t- 1$ others  in  at  least  one  row. We call this  an  {\it $(N,v,t)$-suitable array}. It is clear that $N\geq t$. For example, $\{312645,461523,421365,562134\}$ is $t$-suitable on $[6]$ and its $(4,6,3)$-suitable array is the following.
\[\begin{array}{llllll}
3&1&2&6&4&5\\
4&6&1&5&2&3\\
4&2&1&3&6&5\\
5&6&2&1&3&4\\
\end{array}
\]

The concept of suitable arrays was first introduced by Dushnik \cite{dushnik1950concerning} in 1950 when considering the dimension of partially ordered sets.  Dushnik studied an extremal problem concerning suitable arrays (P1): Given $v$ and $t$, what is the smallest $N$ for which an $(N,v,t)$-suitable array exists? We denote this by $N(v,t)$ \cite{dushnik1950concerning}. Since the $v\times v$ array whose initial elements are $1,2,\ldots,v$ is a $(v,v,t)$-suitable array for each $t\leq v$, we have $N(v,t)\leq v$. So we always assume that $N\leq v$ when we talk about $(N,v,t)$-suitable arrays. By combinatorial arguments, Dushnik \cite{dushnik1950concerning} showed that $N(v,t)=v-j+1$ for each $j$ satisfying $2\leq j\leq \sqrt{v}$ and for each $t$ satisfying \[\left\lfloor\frac{v}{j}\right\rfloor+j-1\leq t<\left\lfloor\frac{v}{j-1}\right\rfloor+j-2.\]This determines $N(v,t)$ exactly for all $t$ in the range $2\lfloor\sqrt{v}\rfloor-1\leq t\leq v$.
%\[\left\lfloor\frac{v}{\lfloor\sqrt{v}\rfloor}\right\rfloor+\lfloor\sqrt{v}\rfloor-1\leq t<v.\] In particular, when the lower bound is attained, both $v$ and $N(v,t)$ grow as $\Theta(t^2)$.

In 1971, Spencer \cite{spencer1971minimal} studied the same problem and showed that for every fixed $t\geq 3$, $N(v,t)\geq \log_2\log_2 v$ and $N(v,t)=O(\log_2\log_2v)$ as $v\rightarrow \infty$. Using probabilistic methods, F\"uredi and Kahn \cite{furedi1986dimensions} showed that  $N(v,t)\leq t^2(1+\log(v/t))$ for all $v$ and $t$ in 1986. Later, Kierstead \cite{kierstead1996order} refined this result in 1996 when $t$ is approximately $\log v$.

In a recent survey paper, Colbourn \cite{colbourn2015suitable} studied the dual extremal problem of suitable sets of permutations (P2): Given $N$ and $t\geq 3$, what is the largest $v$ for which an $(N,v,t)$-suitable array exists? We denote this as SUN$(t,N)$ \cite{colbourn2015suitable}. It is well defined for  $t\geq 3$ by reference to the $(v,v,t)$-suitable array described above,  so SUN$(t,N)\geq N$. Note that SUN$(2,N)$ is not defined since any permutation and its reverse form a $2$-suitable sets for arbitrarily large $v$. By \cite{spencer1971minimal}, we have SUN$(t,N)\leq 2^{2^N}$. Colbourn \cite{colbourn2015suitable} extended this result to SUN$(t,N)=\Theta(2^{2^N})$ for fixed $t$, by linking suitable sets of permutations to  binary covering arrays \cite{lawrence2011survey}. He examined the case when $v$ and $N$ both grow as $t^2$ by making a connection with Golomb rulers and their variants \cite{drakakis2009review,erdos1941problem}. When $t$ is $O(\log N)$, he made a connection with Hadamard matrices \cite{horadam2012hadamard} and Paley matrices \cite{paley1933orthogonal}.

A very interesting observation in Colbourn's paper is that he established an equivalence between a smaller permutation array  and a suitable array. Hence the problem on determining the value SUN$(t,N)$ is transformed by the following  quantity \[\text{SCN}(t,N):=\text{SUN}(t,N)-N,\] whose motivation will be described in Section \ref{sec:pre}.

Colbourn \cite[Section 1]{colbourn2015suitable} derived completely the value of SCN$(t,N)$ when $N< \lfloor\frac{t+1}{2}\rfloor\lceil\frac{t+1}{2}\rceil$, and establish that SCN$(t,\lfloor\frac{t+1}{2}\rfloor\lceil\frac{t+1}{2}\rceil)\geq \lfloor\frac{t}{2}\rfloor+2$. Very recently, Chan and Jedwab \cite{chan2017constructions} proved the other direction of the inequality, which determines SCN$(t,\lfloor\frac{t+1}{2}\rfloor\lceil\frac{t+1}{2}\rceil)= \lfloor\frac{t}{2}\rfloor+2$ for all large $t$. We state this result as follows.

\begin{thm}\label{known1}\cite{chan2017constructions}
\begin{enumerate}[(i)]
\item SCN$(2s+1,(s+1)^2)=s+2$ for all $s\geq 3$.
\item  SCN$(2s,s(s+1))=s+2$ for all $s\geq 2$.
\end{enumerate}
\end{thm}

Motivated by Theorem~\ref{known1}, Chan and Jedwab \cite{chan2017constructions} further considered a question:  whether one can increase the maximum possible value of $v$ from $\lfloor\frac{t}{2}\rfloor+2$  by incrementing the value of the parameter $N=\lfloor\frac{t+1}{2}\rfloor\lceil\frac{t+1}{2}\rceil$ by $1$; or in other words, is  SCN$(t,\lfloor\frac{t+1}{2}\rfloor\lceil\frac{t+1}{2}\rceil+1)> \lfloor\frac{t}{2}\rfloor+2$ for infinitely many $t$? Small examples support this question, for example SCN$(7,17)\geq 6$ and SCN$(9,26)\geq 7$. But surprisingly, the authors in \cite{chan2017constructions} showed that, the value of $N$ can be increased any fixed amount and yet $v$ can be increased from $\lfloor\frac{t}{2}\rfloor+2$ for only finitely many $t$. We state this result as follows.

\begin{thm}\cite{chan2017constructions}\label{known2}
\begin{enumerate}[(i)]
\item For each nonnegative integer $l$, there exists $s_0=s_0(l)$ such that SCN$(2s+1,(s+1)^2+l)=s+2$ for all $s\geq s_0$.
\item For each nonnegative integer $l$, there exists $s_0=s_0(l)$ such that SCN$(2s,$ $s(s+1)+l)=s+2$ for all $s\geq s_0$.
\end{enumerate}
\end{thm}

Theorem~\ref{known2} (i) shows that
\begin{equation}\label{eq1}
\textmd{SCN}(2s+1,(s+1)^2+l)>s+2
\end{equation}
 holds for only finitely many $s$ when $l$ is a fixed positive integer. But if $l$ is allowed to increase with $s$, then (\ref{eq1}) can holds for infinitely many $s$: substitute $s+1$ for $s$ in Theorem~\ref{known1} (ii) to get SCN$(2s+2,(s+1)(s+2))=s+3$, and use that fact that any $(N,v,t)$-suitable core is also an   $(N,v,t-1)$-suitable core (which will be obvious after we give the definition in Section~\ref{sec:pre}), we have SCN$(2s+1,(s+1)(s+2))\geq s+3$, which means $l=s+1$ suffices. Chan and Jedwab \cite{chan2017constructions} then proposed the following problem.

 \vspace{0.5cm}
 {\it Q1: Does there exists a function $l(s)$  which is growing more slowly than linearly with $s$, such that SCN$(2s+1,(s+1)^2+l)\geq s+3$ or SCN$(2s,s(s+1)+l)\geq s+3$ for sufficiently large $s$?}
 \vspace{0.5cm}

 We focus on the problem Q1 in this paper. In fact, we study a problem in a more general pattern.

 \vspace{0.5cm}
 {\it Q2: Let $v=\lfloor\frac{t}{2}\rfloor+\alpha$, where $\alpha\geq 3$ be a fixed constant. Does there exists a function $l(t)$  which is growing more slowly than linearly with $t$, such that SCN$(t,v(t+1-v)+l)\geq v$ for sufficiently large $t$?}
 \vspace{0.5cm}

This paper is organized as follows. In Section~\ref{sec:pre}, we recall the concept of suitable cores as equivalent objects of suitable arrays and some preliminary results.  In Section~\ref{sec:pac}, we give direct constructions of $(N,v,t)$-suitable cores with $v=\lfloor\frac{t}{2}\rfloor+3$, which imply that $l=\Omega(s^{1/3})$ is an answer of Q1 for both cases.  We further generalize our construction by using packings of triples by subsets, and then show that $l=\Omega(t^{1/3})$ is in fact an answer of Q2 for any constant $\alpha\geq 3$. However, this is not the best answer for either Q1 or Q2. In Section~\ref{sec:ram}, we give another construction of $(N,s+3,2s+1)$-suitable cores by applying Ramsey theory, which yields that $l=\Omega(\ln s)$ is an answer of Q1 for this case. We state our first result as follows.
 \begin{thm}\label{upbound}   For each constant $\tau\geq \frac{4}{\ln 2}$,  there exists $s_0=s_0(\tau)$ such that an $(N,s+3,2s+1)$-suitable core exists with $N=(s+3)(s-1)+\tau\ln s$ for all $s\geq s_0$.
\end{thm}

 By Ramsey theory, we also give a nonexistence result of $(N,s+3,2s+1)$-suitable cores when $l=O(\ln s)$, which improves Theorem~\ref{known2} as follows.

\begin{thm}\label{asym}
\begin{enumerate}[(i)]
\item For any function $l=l(s)\leq\frac{\ln s}{6\ln 3}$, SCN$(2s+1,(s+1)^2+l)=s+2$ for all sufficiently large $s$.
\item For any function $l=l(s)\leq \frac{\ln s}{6\ln 3}$,  SCN$(2s,s(s+1)+l)=s+2$ for all sufficiently large $s$.
\end{enumerate}
\end{thm}

 In Section~\ref{sec:gram}, we generalize our construction in Section~\ref{sec:ram} by introducing an extended Ramsey coloring (where each edge is colored by a set of colors), from which we get
 \begin{thm}\label{upbound2}  Let $t=2s+\delta$ and $v=s+\alpha$, where $\delta=0$ or $1$, and $\alpha\geq 3$ is  a fixed constant.
 For each constant $\tau\geq \frac{2r}{\ln r-\ln(r-2)}$ with $r=2\alpha-\delta-2$,  there exists an $(N,v,t)$-suitable core with $N=v(t+1-v)+\tau\ln s$ for all large $s$.
\end{thm}
% We define $N_c(v,t)$ be the smallest $N$ such that there exists an $(N,v,t)$-suitable core.
Theorem~\ref{upbound2} implies that $l=\Omega(\ln s)$ is also an answer of Q2 in general.

%\begin{thm}\label{asym2} Let $\alpha\geq 3$ be a fixed integer and $r=2\alpha-3$.
%\begin{enumerate}[(i)]
%\item For any function $l=l(s)\leq\frac{\ln (s)}{2r\ln r}$, SCN$(2s+1,(s+\alpha)(s+2-\alpha)+l)=s+\alpha$ for all sufficiently large $s$.
%\item For any function $l=l(s)\leq \frac{\ln (s)}{2r\ln r}$,  SCN$(2s,(s+\alpha)(s+1-\alpha)+l)=s+\alpha$ for all sufficiently large $s$.
%\end{enumerate}
%\end{thm}

 \section{Suitable cores}\label{sec:pre}
 This section serves to recast the problem of constructing suitable arrays as the equivalent problem of constructing ``suitable cores'', which is defined by Colbourn \cite{colbourn2015suitable}, who attributed it to Dushnik \cite{dushnik1950concerning}.

 If $A$ is an $(N,v,t)$-suitable array, and $\alpha$ is the initial (or called {\it leader}) element of some row of $A$, then by moving all occurrence of $\alpha$ in all other rows of $A$ to the rightmost positions results in another $(N,v,t)$-suitable array. Applying this to any $(N,v,t)$-suitable array with $N\leq v$ and all leaders in this array, we can get an $(N,v,t)$-suitable array with $N$ different leaders, and all these $N$ elements move to the rightmost $N-1$ positions in the rows starting with  different leaders.

For example, we transform the $(4,6,3)$-suitable array over $[6]$ in Section~\ref{se:intro} to the one on the left. On the right, we have renamed symbols so that the leaders are $3,4,5$ and $6$.

\[\begin{array}{lllllllllllllllllll|ll|llll}
3&1&6&4&2&5&&&&&&&&&&&&&3&1&2&4&6&5\\
4&6&1&3&2&5&&&&&&&&&&&&&4&2&1&3&6&5\\
2&1&6&3&4&5&&&&&&&&&&&&&6&1&2&3&4&5\\
5&6&1&3&4&2&&&&&&&&&&&&&5&2&1&3&4&6\\
\end{array}
\]

Note that the order of the $N$ leaders in the rightmost $N-1$ positions does not affect the suitable property. So to find suitable sets of permutations, it suffices to consider the permutations with the leaders removed. The resulting sets of permutations over $[v-N]$ is the so called  {\it suitable core} \cite{colbourn2015suitable}.

A collection of $N$ permutations over $[v-N]$ is a {\it $t$-suitable core} if it can be extended to an $(N,v,t)$-suitable array by choosing $N$ new symbols, prepending a different one to each permutation, and appending the remaining $N-1$ new symbols in arbitrary order. We denote it by {\it $(N,v-N,t)$-suitable core}. In the example above, the subarray on the right is a $(4,2,3)$-suitable core.

We see in this way that the existence of an $(N,v+N,t)$-suitable array is equivalent to the existence of an $(N,v,t)$-suitable core. Given $N$ and $t$, define SCN$(t,N)$ to be the largest $v$ for which an $(N,v,t)$-suitable core exists. Then SUN$(t,N)=$SCN$(t,N)+N$ provided that $N\geq t$.

To characterize the structure of an $(N,v,t)$-suitable core, we need the following notation.  For an array $C$, symbol $\sigma$ and subset $T$ of other symbols, denote by $C_{pre}(\sigma,T)$ the set of rows of $C$ for which $\sigma$ either starts a row or is preceded only by elements of $T$. In other words, $C_{pre}(\sigma,T)$  is the set of rows of $C$ where $\sigma$ precedes all elements of $[v]\setminus(T\cup\{v\})$.

 \begin{pro}\label{suitablecore}\cite{colbourn2015suitable,chan2017constructions}
 Let $C$ be an $N\times v$ array. Then the following statements are equivalent:
 \begin{enumerate}[(i)]
\item $C$ is an $(N,v,t)$-suitable core.
\item For each $s$ satisfying $0\leq s\leq t-1$, each symbol of $C$ precedes each subset of $s$ others in at least $t-s$ rows.
\item For each symbol $\sigma$ of $C$ and for each subset $T$ of other symbols, $|C_{pre}(\sigma,T)|$ $\geq t+1-v+|T|$.
\end{enumerate}
 \end{pro}

Besides Proposition~\ref{suitablecore}, the following lemma is very useful in the existence and non-existence proofs of suitable cores.

 \begin{lem}\label{numrows}\cite{chan2017constructions}
Suppose that $C$ is an $(N,v,t)$-suitable core.
 \begin{enumerate}[(i)]
\item Let $v\leq t$. Then each $k\in[v]$ starts a row at least $t+1-v$ times.
\item Let $v\leq t+1$, let  $j$ and $k$ be two different symbols each  starting exactly $t+1-v$ rows. Then there is at least one row that starts with $jk$.
    \item Let $v\leq t+2$, let $k$ be a symbol which starts exactly $t+2-v$ rows, and let $i,j$ be two other distinct symbols. If neither $ik$ nor $jk$ starts a row, then there is at least one row that starts with $ijk$ or $jik$.
\end{enumerate}
 \end{lem}

From Lemma~\ref{numrows}, the following result is obvious and will be used repeatedly in our constructions and proofs.
 \begin{lem}\label{compl}
Suppose that $C$ is an $(N,v,t)$-suitable core over $[v]$. Let $R$ be the set of elements starting a row more than $t+1-v$ times. For each $i\in [v]$, let $B_i$ be the collection of elements $j$ such that $ij$ does not start a row in $C$. Then $B_i\subset R$ for each $i\in [v]$, and $|B_i|\geq 2v-t-2$ for each $i\in [v]\setminus R$.
 \end{lem}
\begin{pf}By Lemma~\ref{numrows} (i), each $j\in [v]\setminus R$  starts a row exactly $t+1-v$ times.  For each such $j\in [v]\setminus R$, and for each different $i\in [v]$ , $ij$ must start a row at least once by Lemma~\ref{numrows} (ii). Hence, $B_i\subset R$ for each $i\in [v]$. It is obvious that $|B_i|\geq(v-1)-(t+1-v)=2v-t-2$ for each $i\in [v]\setminus R$. \qed
\end{pf}

The following result links suitable cores with parameters $t$ and $t+1$.

 \begin{lem}\label{trans}\cite{chan2017constructions}
Suppose that $SCN(t,N)\geq v$ and $N>v(t+1-v)$. Then $SCN(t+1,N+v-1)\geq v$.
 \end{lem}

 Since each $k\in[v]$ starts a row at least $t+1-v$ times, we have $N\geq v(t+1-v)$ in an $(N,v,t)$-suitable core. In the remaining of this paper, we always write $N$ in the form of $v(t+1-v)+l$. The main problem in our constructions is how to set leaders of the remaining $l$ rows.

\section{Constructions from packings of triples}\label{sec:pac}

In this section, we give explicit constructions of $(N,v,t)$-suitable cores with $(N,v,t)=((s+3)(s-1)+l,s+3,2s+1)$ and $((s+3)(s-2)+l,s+3,2s)$, where $l=\Omega(s^{1/3})$. This affirms the question Q1 proposed by Chan and  Jedwab. Further, we extend our constructions to  $(N,v,t)$-suitable cores with $(N,v,t)=((s+\alpha)(s+\delta-\alpha+1)+\Omega(s^{1/3}),s+\alpha,2s+\delta)$ for all absolute constants  $\alpha\geq 3$ and $\delta=0,1$. This gives an answer of Q2.

\subsection{Case 1: $t$ is odd}

Let $t=2s+1$ and $v=s+3$, we will construct an $(N,v,t)$-suitable core with $N=(s+3)(s-1)+l$ such that ${l \choose 3}\geq s+3\geq l$.

We first briefly describe our main idea of the construction. By Lemma~\ref{numrows} (i), each of the $s+3$ symbols of $C$ starts a row at least $t+1-v=s-1$ times. This accounts for $(s+3)(s-1)$ rows of $C$, leaving $l$ rows to account for. Now we let the $l$ remaining rows start with different symbols. In other words, we will construct an $N\times v$ array with each symbol from $[v-l]$ starting a row exactly $s-1$ times, and each symbol from $[v-l+1,v]$ starting a row exactly $s$ times. Let $c=v-l$ and $R=[c+1,v]$, and hence $|R|=l$. We use the same notation $B_i$ as in Lemma~\ref{compl}.  Then $B_i\subset R$ for each $i\in [v]$. Further, we assume that each pair $ij$ starts a row at most once. Thus $|B_i|=v-1-(s-1)=3$ for each $i\in [c]$ and $|B_i|=v-1-s=2$ for each $i\in R$. In our construction, we will use certain subsets $B_i\subset R$, $i\in [v]$ to define the first two elements of each row of $C$, such that the partial array can be completed to an $(N,v,t)$-suitable core.

\begin{cons}\label{con1} Let $t=2s+1$, $v=s+3$ and $N=(s+3)(s-1)+l$ satisfying that  ${l \choose 3}\geq s+3\geq l$. Let $c=v-l$ and $R=[c+1,v]$. Let $B'_i$, $i\in [v]$ be a set of $v$ distinct $3$-subsets of $R$, such that $B'_i$ contains the symbol $i$ for each $i\in R$.  Then let $B_i=B'_i$ if $i\in [c]$ and $B_i=B_i'\setminus\{i\}$ if $i\in R$. We construct an $N\times v$ array $C$ as follows. For each $i\in [c]$, and each $j\not\in B_i$, $ij$ starts a row of $C$ exactly once. Next we assign the third elements for some rows of $C$ as follows, and then complete each row arbitrarily to a permutation over $[v]$.
 \begin{enumerate}[(O1)]
\item  If $B_i\cap B_j=\{k_1,k_2\}$, then let $k_1$ appear third after $ij$ and $k_2$ appear third after $ji$ in $C$.
\item If $B_i\cap B_j=\{k\}$, then let $k$ appear third after  $ij$ or $ji$ in $C$.
\end{enumerate}
\end{cons}

The following lemma shows that Construction~\ref{con1} gives a suitable core.

\begin{lem}\label{triples} There exists an  $(N,s+3,2s+1)$-suitable core with $N=(s+3)(s-1)+l$ provided that ${l \choose 3}\geq s+3\geq l$, for all $s\geq 1$.
\end{lem}
\begin{pf} We first prove that the two operations (O1) and (O2) in Construction~\ref{con1} are always executable.

If $B_i\cap B_j=\{k_1,k_2\}$, we need to show that both $ij$ and $ji$ start a row in $C$.  It is true if $i,j$ are both in $[c]$, since $B_i,B_j\subset R$. If $i\in [c]$ and $j\in R$, then $ji$ starts a row trivially. Further, $B_j'=\{j,k_1,k_2\}$, and hence $j\not \in B_i$ since otherwise $B_i'=B_j'=\{j,k_1,k_2\}$. This implies that $ij$ starts a row.  If both $i,j\in R$, then $B_i=B_j=\{k_1,k_2\}$, that is, $i\not\in B_j$ and $j\not\in B_i$ implying that both $ij$ and $ji$ start a row.

If $B_i\cap B_j=\{k\}$, we need to show that at least one of the pairs $ij$ and $ji$ starts  a row in $C$. It is true if one of the symbols $i,j$ is in $[c]$.  If both $i,j\in R$, then $i\in B_j$ and $j\in B_i$ can not happen simultaneously, since otherwise $B_i'=B_j'=\{i,j,k\}$, a contradiction.

Now we prove that the $N\times v$ array $C$ is actually an $(N,s+3,2s+1)$-suitable core using Proposition~\ref{suitablecore}. Let  $\sigma\in [v]$ and let $T$ be a (possibly empty) set of symbols other than $\sigma$. We distinguish two cases.

 \begin{enumerate}[C1]
\item  $\sigma\in [c]$. Then $\sigma$ starts a row exactly $v-1-|B_{\sigma}|=v-1-3=t+1-v$ times. Since $B_j\subset R$ for each $j\in T$, $j\sigma$ starts a row in $C$. So we have $|C_{pre}(\sigma,T)|\geq t+1-v+|T|$.
    \item $\sigma\in R$. Then $\sigma$ starts a row exactly $v-1-|B_{\sigma}|=v-1-2=t+2-v$ times. If $|T|\leq 1$, then $|C_{pre}(\sigma,T)|\geq t+2-v\geq t+1-v+|T|$. Assume that $T=\{a_1,\ldots,a_g\}$ with $g\geq 2$. Let $T'\subset T$ be the collection of symbols $a_i$ such that $\sigma\in  B_{a_i}$. For each pair $\{a_i,a_j\}\subset T'$,  we have either $a_ia_j\sigma$ or $a_j a_i\sigma$ starts a row by Construction~\ref{con1}. Hence $|C_{pre}(\sigma,T)|\geq t+2-v+|T\setminus T'|+{|T'|\choose 2}\geq t+1-v+|T|$.\qed
\end{enumerate}

\end{pf}

\subsection{Case 2: $t$ is even}
Let $t=2s$, $v=s+3$ and $N=(s+3)(s-2)+l$ whose conditions will be given later. Similar to Construction~\ref{con1}, we will use a set of quadruples to define an $(N,v,t)$-suitable core, where these quadruples form a $3$-$(l, 4, 1)$ packing. We first introduce the concept of packings from combinatorial design theory \cite{douglas2007packings}.

Let $l \geq k \geq  t$ and $\lambda$ be positive integers. A $t$-$(l, k, \lambda)$ packing is a pair $(X, \LB)$, where $X$ is an $l$-set of elements
({\it points}) and $\LB$ is a collection of $k$-subsets of $X$ ({\it blocks}), such that every $t$-subset
of points occurs in at most $\lambda$ blocks in $\LB$.
Given $t$, $k$, and $l$, the determination of the packing number $D(l, k, t)$, the maximum size of a $t$-$(l, k, 1)$ packing, constitutes a central problem
in combinatorial design theory, as well as in coding theory \cite{mills1992coverings}.
When $k=4$ and $t=3$, the value of  $D(l, 4, 3)$ has been completely determined by constructive methods, see \cite{Hanani1960CJM,hanani1979class,Ji2006DCC,BJ2014}, and it achieves the well known Johnson bound given below:

\begin{equation}\label{eq2} D(l,4,3)=\begin{dcases}
\left\lfloor\frac{l}{4}\left\lfloor\frac{l-1}{3}\left\lfloor\frac{l-2}{2}\right\rfloor\right\rfloor\right\rfloor, &\text{if $l\not\equiv 0 \pmod 6$;} \\
\left\lfloor\frac{l}{4}\left(\left\lfloor\frac{l-1}{3}\left\lfloor\frac{l-2}{2}\right\rfloor\right\rfloor-1\right)\right\rfloor, &\text{if $l\equiv 0 \pmod 6$. }
\end{dcases}\end{equation}

\begin{cons}\label{con2} Let $t=2s$, $v=s+3$ and $N=(s+3)(s-2)+l$ satisfying that  $D(l, 4, 3)\geq s+3\geq l$. Let $c=v-l$ and $R=[c+1,v]$. Let $B'_i$, $i\in [v]$ be blocks of a $3$-$(l, 4, 1)$ packing over $R$, such that $B'_i$ contains the symbol $i$ for each $i\in R$. This can be done  when $s\geq 4$. Let $B_i=B'_i$ if $i\in [c]$ and $B_i=B_i'\setminus\{i\}$ if $i\in R$. Then we construct an $N\times v$ array $C$ by these sets $B_i$, $i\in [v]$, using the same method in Construction~\ref{con1}.
\end{cons}

Construction~\ref{con2} produces an  $(N,s+3,2s)$-suitable core with the prescribed parameters.

\begin{lem}\label{quadruples} There exists an  $(N,s+3,2s)$-suitable core with $N=(s+3)(s-2)+l$ provided that $D(l, 4, 3)\geq s+3\geq l$, for all $s\geq 4$.
\end{lem}
\begin{pf} Note that $|B_i\cap B_j|\leq |B'_i\cap B'_j|\leq 2$ for any two distinct $i,j$ by the definition of $3$-$(l, 4, 1)$ packings. Similar to Lemma~\ref{triples}, we need to prove that the two operations (O1) and (O2) in Construction~\ref{con1} are always executable, under the assumption given in Construction~\ref{con2}.

If $B_i\cap B_j=\{k_1,k_2\}$, we need to show that both $ij$ and $ji$ start a row in $C$. Since $B_i,B_j\subset R$, it is true if $i,j$ are both in $[c]$. If $i\in [c]$ and $j\in R$, then  $ji$ starts a row since $i\not\in B_j$. Further, $j\not \in B_i$, since otherwise $B_i'\cap B_j'=\{j,k_1,k_2\}$. This implies that $ij$ starts a row.  Now assume that $i,j\in R$. If  $i\in B_j$ then $B_i'\cap B_j'=\{i,k_1,k_2\}$; if $j\in B_i$, then $B_i'\cap B_j'=\{j,k_1,k_2\}$,  contradictions. So $i\not\in B_j$ and $j\not\in B_i$, which means  that both $ij$ and $ji$ start a row.

If $B_i\cap B_j=\{k\}$, we need to show that at least one of the pairs $ij$ and $ji$ starts a row in $C$. It is true if one of the symbols $i,j$ is in $[c]$.  If $i,j\in R$, then $i\in B_j$ and $j\in B_i$ can not happen simultaneously, since otherwise $B_i'\cap B_j'=\{i,j,k\}$, a contradiction.

The proof that the $N\times v$ array $C$ is actually an $(N,s+3,2s)$-suitable core is similar to that in Lemma~\ref{triples}, thus is omitted.\qed

\end{pf}

\subsection{General case: $v=\lfloor\frac{t}{2}\rfloor+\alpha$ with $\alpha\geq 3$}

In this subsection, we generalize the methods in Constructions~\ref{con1} and~\ref{con2} to the case that $t=2s+0$ or $1$, and $v=s+\alpha$ with $\alpha\geq 3$ being an absolute constant.
We will construct an $(N,v,t)$-suitable core with $N=v(t+1-v)+\Omega(s^{1/3})$.

\begin{cons}\label{con3} Given integers $\delta=0$ or $1$, and $\alpha\geq 3$, let $t=2s+\delta$, $v=s+\alpha$ and $N=v(t+1-v)+l$ satisfying that  $D(l, k, 3)\geq v\geq l$, where $k=2\alpha-\delta-2$. Let $c=v-l$ and $R=[c+1,v]$. Let $B'_i$, $i\in [v]$ be blocks of a $3$-$(l, k, 1)$ packing over $R$, such that $B'_i$ contains the symbol $i$ for each $i\in R$. This can be done  when $s$ is large. Then let $B_i=B'_i$ if $i\in [c]$ and $B_i=B_i'\setminus\{i\}$ if $i\in R$. Then we construct an $N\times v$ array $C$ by these sets $B_i$, $i\in [v]$, using the same method in Construction~\ref{con1}.
\end{cons}

The proof of the following result is similar to that of Lemma~\ref{quadruples}, which we leave to readers.

\begin{lem}\label{ksubsets} Let $\delta=0$ or $1$, and $\alpha\geq 3$ be integers. Then for large $s$, there exists an  $(N,s+\alpha,2s+\delta)$-suitable core with $N=(s+\alpha)(s+\delta-\alpha+1)+l$ provided that $D(l, 2\alpha-\delta-2, 3)\geq s+\alpha\geq l$.
\end{lem}

Note that the packing number $D(l, k, 3) \leq {l\choose 3}/{k\choose 3}$. R\"{o}dl \cite{rodl1985packing} was
the first to show that this upper bound can be attained asymptotically. %Let $\epsilon_{t,k}(l)$  be the fraction of $t$-subsets not contained
%in any blocks of a $t$-$(l, k, 1)$ packing of maximum size. In other words, $D(l, k, t)= (1 - \epsilon_{t,k}(l)){l\choose t}/{k\choose t}$.
That is, when $l$ is large enough, we have a $3$-$(l, k, 1)$ packing with number of blocks arbitrarily close to the upper bound. Hence, we can conclude from  Lemma~\ref{ksubsets}, that an  $(N,s+\alpha,2s+\delta)$-suitable core exists with $N=(s+\alpha)(s+\delta-\alpha+1)+\Omega(s^{1/3})$ when $s$ is sufficiently large.

\section{Proofs of Theorems~\ref{upbound} and~\ref{asym}}\label{sec:ram}
In this section, we apply Ramsey theory to prove existence and nonexistence results of suitable cores. Similar ideas have been used in \cite{chan2017constructions} to prove the nonexistence part of Theorem~\ref{known2}.
%
%Ramsey Theory studies conditions when a combinatorial object contains necessarily some
%smaller given objects. The role of Ramsey numbers is to quantify some of the general existen-
%tial theorems in Ramsey Theory.

We will make use of the following notation. Let $G$ be a graph, $V(G)$ the set of vertices
of $G$, and $E(G)$ the set of edges of $G$. An $r$-coloring, $\chi$, will be assumed to be an edgewise
coloring, i.e. $\chi(G)$ : $E(G)\rightarrow \{1, 2,\ldots,r\}$.  We denote by $K_n$ the complete graph on $n$ vertices.

 Let $r \geq 2$, and let $k_i \geq 2$, $1 \leq i \leq r$. The {\it Ramsey  number} $R(k_1, k_2,\ldots,k_r)$ is
defined to be the minimal integer $n$ such that any edgewise $r$-coloring of $K_n$ must contain,
for some $j$, $1 \leq j \leq r$, a monochromatic $K_{k_j}$ of color $j$ (that is, all edges in this clique have color $j$). If we are considering the
{\it diagonal} Ramsey numbers, i.e. $k_1=k_2=\cdots=k_r=k$, we will use $R(k;r)$ to denote the
corresponding Ramsey number. A {\it Ramsey $r$-coloring} for $R(k_1, k_2,\ldots,k_r)$ is an $r$-coloring of the
complete graph on $n <R(k_1, k_2,\ldots,k_r)$ vertices which does not admit any monochromatic $K_{k_j}$ subgraph
of color $j$ for $j = 1, 2,\ldots,r$. Note that a  Ramsey $r$-coloring for $R(k_1, k_2,\ldots,k_r)$ of $K_n$ exists if and only if $n < R(k_1, k_2,\ldots,k_r)$.

 In the case of two colors ($r = 2$) one deals with classical graph Ramsey numbers, which have
been studied extensively for 50 years. Much less has been done for multicolor numbers
($r \geq 3$). The significant lower bound for the diagonal Ramsey number $R(k,k)\geq \frac{1}{e\sqrt{2}}k2^{\frac{k}{2}}$
was proved by Erd\"{o}s \cite{erdos1947some} in 1947 by using probabilistic method.  An easy extension of the Erd\"{o}s-Szekeres
argument \cite{erdios1935combinatorial} gives an upper bound for the multicolour diagonal Ramsey
number of the form $R(k; r) \leq r^{ rk}$, see \cite{greenwood1955combinatorial}. We include a lower bound recurrence  found by Robertson in \cite{robertson1999new} and \cite{robertson2002new}: for $k, l \geq 3$, we have $R(3,k,l)\geq 4R(k,l-2)-3$.

For more results on the known bounds on various types of Ramsey numbers, see \cite{conlon2015recent} by Conlon et al. and a regularly updated survey \cite{radziszowski1994small} maintained by Radziszowski.
%\begin{equation}\label{eq3}
%R(k_1,k_2,\ldots,k_r)\geq (k_1-1)(R(k_2,\ldots,k_{r})-1)
%\end{equation}

%
%Stanis{\l}aw Radziszowski maintains a\cite{radziszowski1994small,conlon2015recent} of the most
%recent results on the best known bounds on various types of Ramsey numbers.

The construction we give below is quite different from the ones given in Section~\ref{sec:pac}. In Section~\ref{sec:pac}, we construct  ``balanced'' suitable cores, that is,   the number of rows starting with different symbols are  almost equal. But in the suitable cores constructed below,  most symbols start a row with the least necessary number of times, and the remaining symbols start a row far more times.

\begin{cons}\label{con4} Let $t=2s+1$, $v=s+3$ and $N=(s+3)(s-1)+l$.   We construct an $N\times v$ array $C$ as follows.  For each $i\in [s]$ and each $j\in [v]$ with $j\neq i$, $ji$ starts a row. This accounts for $(s+3)(s-1)+3$ rows. Let $R=\{s+1,s+2,s+3\}$. We assume that symbol $s+1$ starts $k_1-1$ other rows, symbol $s+2$ starts $k_2-1$ other rows and symbol $s+3$ starts $k_3-1$ other rows, such that $k_1,k_2,k_3\geq 3$, and $k_1+k_2+k_3=l$. Now we have in total $N$ rows. We further assume that for each $i\in R$, and each different $j\in [v]$, $ij$ starts a row at least once.

 Construct a complete graph $G=K_s$ with vertex set $[s]$. Suppose that  $s< R(k_1+1,k_2+1,k_3+1)$, then we have a  Ramsey $3$-coloring for $R(k_1+1,k_2+1,k_3+1)$  of $K_s$, that is, there does not exist any monochromatic $K_{k_h+1}$ subgraph of color $h$, for $h = 1, 2,3$.  Now we assign the third elements for the  $s(s-1)$ rows of $C$ starting by $ij$, with $i,j\in [s]$. For each pair $i,j\in [s]$, there are exactly two rows starting with $ij$ or $ji$. If the edge $i,j$ in $G$ is colored by $h$, then assign the third elements of these two row by the two elements $s+h'$, with $h'\in [3]\setminus \{h\}$.  Finally, complete each row arbitrarily to a permutation over $[v]$.
%
%  by $s+1$, $s+2$ and $s+3$ as follows, and then complete each row arbitrarily to a permutation over $[v]$.
%Suppose then we can assign the third element of the first $s(s-1)$ rows as follows.
%
% For each pair $i,j\in [s]$, color the edge between $i,j$ with color $h$ if neither $ij(s+h)$ nor $ji(s+h)$ starts a row, where $h\in [3]$.
%
% For each pair of distinct elements $i, j\in [s]$, we assume that two different symbols from $s+1$, $s+2$ and $s+3$  follow $ij$ or $ji$. That is, for each $i,j$, exactly one of $s+1$, $s+2$ and $s+3$ does not come immediately after the rows starting with $ij$ or $ji$.
\end{cons}

We show in the following lemma that Construction~\ref{con4} gives a suitable core.

\begin{lem}\label{ramsey1} Let $l$ be a positive integer such that there exists three integers $k_1,k_2,k_3\geq 3$, and $k_1+k_2+k_3=l$. Then there exists an  $(N,s+3,2s+1)$-suitable core with $N=(s+3)(s-1)+l$ provided that $s< R(k_1+1,k_2+1,k_3+1)$.
\end{lem}
\begin{pf} We again prove it by using Proposition~\ref{suitablecore}. Let  $\sigma\in [v]$ and let $T$ be a  set of symbols other than $\sigma$. We distinguish two cases.

 \begin{enumerate}[C1]
\item  $\sigma\in [s]$. Since for each $j\in T$, $j\sigma$ starts a row, we have $|C_{pre}(\sigma,T)|\geq s-1+|T|=t+1-v+|T|$.
    \item $\sigma=s+h$, for some $h\in [3]$. Then $\sigma$ starts a row exactly $s+k_h-1=t+1-v+k_h$ times. If $|T|\leq k_h$, then $|C_{pre}(\sigma,T)|\geq t+1-v+k_h\geq t+1-v+|T|$. Let $T=\{a_1,\ldots,a_g\}$ with $g\geq k_h+1$. Let $T'\subset T$ be the collection of symbols $a_i$ such that $a_i\sigma$ does not start a row. Hence $T'\subset [s]$, since for each $i\in R$, $i\sigma$
starts a row at least once. Now for each $(k_h+1)$-subset $\{a_{i_1},\ldots, a_{i_{k_h+1}}\}\subset T'$,  we have at least one pair $a_{i_x},a_{i_y}$ such that the edge between them does not have color $h$ in $G$. The reason here is that there is no monochromatic $K_{k_h+1}$ in color $h$ in the Ramsey $3$-coloring of $K_s$. Hence, either $a_{i_x}a_{i_y}\sigma$ or $a_{i_y}a_{i_s}\sigma$ starts a row in $C$ by Construction~\ref{con4}. This means that each  $(k_h+1)$-subset of $T'$ contributes $|C_{pre}(\sigma,T)|$ at least one.
         Hence $|C_{pre}(\sigma,T)|\geq t+1-v+k_h+|T\setminus T'|+{|T'|\choose {k_h+1}}$. We can split into two cases: $|T'|\leq k_h$ or $|T'|\geq k_h+1$. But for both cases, we have $|C_{pre}(\sigma,T)| \geq t+1-v+k_h+|T|-|T'|+{|T'|\choose {k_h+1}}\geq t+1-v+|T|$. \qed
\end{enumerate}
\end{pf}
%We introduce the {\it Lambert W function},  which is a set of functions, namely the branches of the inverse relation of the function $f(z) = ze^z$,  and $z$ is any complex number. In other words, if $z=xe^x$, then $x=W(z)$. Here, we focus on the positive real-valued $W$, which is increasing by $z$. In this case, $W(z)<\ln(z)$. It can be extended to the function $z=xa^x$ using the identity $x=\frac{W(\ln(a)z)}{\ln(a)}$. It is shown \cite{} that the following bound holds for $x\geq e$: \[\ln x -\ln \ln x + \frac{\ln \ln x}{2\ln x} \le W_0(x) \le \ln x - \ln\ln x + \frac{e}{e - 1} \frac{\ln \ln x}{\ln x.}\]

Now we prove Theorem~\ref{upbound}. For  simplicity, we assume that all the defined parameters are integers in the remaining of this paper.

%\begin{thm}\label{upbound}   For each constant $\sigma\geq \frac{4}{\ln 2}$,  there exists $s_0=s_0 (\delta)$ such that an $(N,s+3,2s+1)$-suitable core with $N=(s+3)(s-1)+\sigma\ln (s)$ exists for all $s\geq s_0$.
%\end{thm}
\begin{pot3}Let $l=\tau\ln s$ and set $k_1=3,k_2=(l-5)/2$ and $k_3=(l-1)/2$ in Lemma~\ref{ramsey1}. Let $k=(l-3)/2$. Then an $(N,s+3,2s+1)$-suitable core with $N=(s+3)(s-1)+l$ exists if $s< R(4,k,k+2)$ by Lemma~\ref{ramsey1}. Note that $R(4,k,k+2)\geq R(3,k,k+2)\geq 4R(k,k)-3\geq 4\cdot \frac{1}{e\sqrt{2}}k2^{\frac{k}{2}}-3\geq  \frac{k}{2}2^{\frac{k}{2}}\geq 2^{\frac{k}{2}}$.  When $s\geq s_0=2^{3/(\tau\ln2-4)}$, we have  $s\leq 2^{\frac{k}{2}}< R(4,k,k+2)$.\qed
\end{pot3}

Next, we show the other side of Theorem~\ref{upbound}, that is, when $\tau$ is small enough, there exists an  $(N,s+3,2s+1)$-suitable core with $N=(s+3)(s-1)+\tau\ln s$ only for finitely many $s$. We need the following lemma from \cite{chan2017constructions}.

\begin{lem}\label{midl}\cite{chan2017constructions}
Let $e$ and $m$ be positive integers and $d\geq (e-1)(2m+1)+1$. Let $A$ be a set of size $d$ and let $g$ be a function from $A$ to subsets of $A$ of size at most $m$. Then there exists a subset $B$ of $A$ of size $e$ such that $j\not \in g(i)$ for all distinct $i,j\in B$.
\end{lem}
\begin{lem}\label{lbound}   There exists $s_0$ such that an  $(N,s+3,2s+1)$-suitable core with $N=(s+3)(s-1)+\frac{\ln s}{6\ln 3}$ does not exist for all $s\geq s_0$.
\end{lem}
\begin{pf}  Let $l= \frac{\ln s}{6\ln 3}$, $v=s+3$ and $t=2s+1$. Suppose, for a contradiction, that there exists an $(N,s+3,2s+1)$-suitable core $C$ with $N=(s+3)(s-1)+l$.

Without loss of generality, we assume that $C$ is over $[v]$ and the number of rows starting with $i$ is nondecreasing with $i$. By Lemma~\ref{numrows} (i), each of the $s+3$ symbols of $C$ starts a row at least $s-1$ times. This accounts for $(s+3)(s-1)$ rows of $C$, leaving $l$ more rows to account for. Let $c$ be the number of symbols starting a row exactly $s-1$ times, that is, each symbol of $[c]$ starts a row exactly $s-1$ times. Then $c\geq v-l$. Let $R=[c+1,v]$, which consists of symbols starting more than $s-1$ rows of $C$. Let $m=|R|= v-c\leq l$.

By Lemma~\ref{numrows} (ii), $C$ contains a row starting with $ij$ for each $j\in[c]$ and for each $i\neq j$. Hence, $c-1\leq s-1$, that is $c\leq s$.  For each $j\in [c]$, let $B_j$ be the set of elements $i$, such that $i$ does not appear second after $j$ in $C$. Then by Lemma~\ref{compl}, for each $j\in [c]$,  $B_j\subset R$, and $|B_j|\geq (v-1)-(s-1)=v-s=3$. Now let $f$ be  a function from $[c]$ to $3$-subsets of $R$, such that $f(j)\subset B_j$, $j\in [c]$.

Let $d=3ls^{1/2}$. Since $c\geq v-l=s+3-l$, we have $c\geq {l\choose 3}(d-1)+1\geq {m\choose 3}(d-1)+1$ when $s$ is large enough.
Then by pigeonhole principle, there exists a set of $d$ numbers $A=\{a_1,a_2,\ldots,a_d\}\subset [c]$ for which $f(a_1)=f(a_2)=\cdots=f(a_d)$. Let $\{r_1,r_2,r_3\}=f(a_1)$.

Next, define a function $g$ from $A$ to subsets of $A$ as follows: for each $a\in A$, let $g(a)$ be a set consisting of elements  of $A$ that occurs second after $a$ at least twice in $C$. So $|g(a)|\leq s-1-(c-1)=s-c\leq l$. Let $e=s^{1/2}$, then $d\geq(e-1)(2l+1)+1$  when $s$ is large enough. By Lemma~\ref{midl}, there exists a subset $B=\{b_1,b_2,\ldots,b_e\}$ of $A$ of size $e$ such that $b_i\not\in g(b_j)$ for all distinct $i,j$. It follows that for each pair of distinct elements $b_i,b_j$ of $B$ there is exactly one row of $C$ starting $b_ib_j$.

Now construct  a complete graph $G$ with vertex set $B$. For each  $b\in B$, observe that $f(b)=\{r_1,r_2,r_3\}$. Then for each $b_i,b_j\in B$, there is at least one element of the set $\{r_1,r_2,r_3\}$ that precedes the other two in neither the row starting $b_ib_j$ nor the row staring $b_jb_i$; choose one such element and color the edge between vertices $b_i,b_j$ with color $1$ if the choice is $r_1$, color $2$ if it is $r_2$, and color $3$ if it is $r_3$. So the resulting graph $G$ is an edgewise $3$-coloring of $K_e$.

Suppose in $C$, the symbols $r_1,r_2,r_3$ start a row $s-1+k_1$, $s-1+k_2$ and $s-1+k_3$ times, respectively, for some positive integers $k_1,k_2,k_3$. Note that $k_1,k_2,k_3\leq l-1$. Since $e\geq 3^{3l}\geq R(l,l,l)\geq R(k_1+1,k_2+1,k_3+1)$, there does not exist a Ramsey $3$-coloring for $R(k_1+1,k_2+1,k_3+1)$ of $K_e$. That is, there must exist a monochromatic $K_{k_i+1}$ in color $i$, for some $i\in [3]$ in $G$. Suppose the monochromatic $K_{k_i+1}$ is with vertex set $T\subset B\subset A\subset [c]$, and $|T|=k_i+1$. Then $|C_{pre}(r_i,T)|=s-1+k_i<t+1-v+|T|$, which is a contradiction by Proposition~\ref{suitablecore} (iii). \qed
\end{pf}

Now we are ready to prove Theorem~\ref{asym}.

%\begin{thm}\label{asym}
%\begin{enumerate}[(i)]
%\item For any function $l=l(s)$ growing less than $\frac{\ln (s)}{6\ln 3}$, SCN$(2s+1,(s+1)^2+l)=s+2$ for all sufficiently large $s$.
%\item For any function $l=l(s)$ growing less than $\frac{\ln (s)}{6\ln 3}$,  SCN$(2s,s(s+1)+l)=s+2$ for all sufficiently large $s$.
%\end{enumerate}
%\end{thm}
\begin{pot4} The statement (i) is a direct consequence of  Theorem~\ref{known2} (i) and Lemma~\ref{lbound}.
For (ii), if SCN$(2s,s(s+1)+l)\geq s+3$, then SCN$(2s+1,(s+1)^2+l+1)\geq s+3$ by Lemma~\ref{trans}, a contradiction by taking $l'=l+1$.\qed
\end{pot4}

\section{Proof of Theorem~\ref{upbound2}}\label{sec:gram}

In order to generalize Construction~\ref{con4}, we need an extension of Ramsey number to multicolors for each edge.

Let $r>m\geq 1$ be  positive integers. An  {\it $(r;m)$-coloring},  $L(G)$ : $E(G)\rightarrow {[r]\choose m}$ is a function assigning to each edge $e\in E(G)$ a list of $m$ colors $L(e)\subset [r]$. Under this definition, if  $i\in L(e)$ for all edges $e$ of a complete subgraph $K_k$ of $G$, then we say that $G$ contains a monochromatic  $K_{k}$ in color $i$. A {\it Ramsey $(k_1, k_2,\ldots, k_r)^m$-coloring} of $K_n$, $k_i \geq  1$, is an $(r;m)$-coloring such that it does not contain any monochromatic complete subgraph
$K_{k_i}$ in color $i$, for $1 \leq i \leq  r$.  For example, the graph in Fig.~\ref{fig1} has two monochromatic $K_3$'s in color $1$ with vertex sets $\{1,2,3\}$ and $\{1,3,4\}$, and the one in Fig.~\ref{fig2} is a Ramsey $(3, 3,3)^2$-coloring of $K_4$.

\begin{figure}[!htb]
    \centering
    \begin{minipage}{.5\textwidth}
        \centering
\begin{tikzpicture}
\coordinate (3) at (0,0);
\coordinate (4) at (3,0);
\coordinate (1) at (0,2);
\coordinate (2) at (2,2);
\filldraw
(3) circle (1pt) node[below] {3} --node[font=\fontsize{8}{0}\selectfont,left=1pt] {$1,2$}
(1) circle (1pt) node[above] {1}     --node[font=\fontsize{8}{0}\selectfont,above] {$1,3$}
(2) circle (1pt) node[above] {2} --node[font=\fontsize{8}{0}\selectfont,right] {$2,3$}
(4) circle (1pt) node[below] {4};
\draw (1)-- node[font=\fontsize{8}{0}\selectfont,right] {$1,2$} (4)-- node[font=\fontsize{8}{0}\selectfont,below] {$1,3$}(3)--node[font=\fontsize{8}{0}\selectfont,below] {$1,2$}(2);
\end{tikzpicture}
     \caption{Non-Ramsey coloring}
        \label{fig1}
  \end{minipage}%
    \begin{minipage}{0.5\textwidth}
        \centering
\begin{tikzpicture}
\coordinate (3) at (0,0);
\coordinate (4) at (3,0);
\coordinate (1) at (0,2);
\coordinate (2) at (2,2);
\filldraw
(3) circle (1pt) node[below] {3} --node[font=\fontsize{8}{0}\selectfont,left=1pt] {$2,3$}
(1) circle (1pt) node[above] {1}     --node[font=\fontsize{8}{0}\selectfont,above] {$1,3$}
(2) circle (1pt) node[above] {2} --node[font=\fontsize{8}{0}\selectfont,right] {$2,3$}
(4) circle (1pt) node[below] {4};
\draw (1)-- node[font=\fontsize{8}{0}\selectfont,right] {$1,2$} (4)-- node[font=\fontsize{8}{0}\selectfont,below] {$1,3$}(3)--node[font=\fontsize{8}{0}\selectfont,below] {$1,2$}(2);
\end{tikzpicture}
     \caption{Ramsey coloring}
        \label{fig2}
    \end{minipage}

\end{figure}

  The extended   Ramsey
number $R^m(k_1, \ldots, k_r)$ is defined to be the least integer $n > 0$ such that there is no Ramsey $(k_1, k_2,\ldots, k_r)^m$-coloring of $K_n$.  When $k_1=k_2=\cdots=k_r=k$, we simply write $R^m(k;r)$. The case of $m=1$ is the classical Ramsey numbers with multicolors.

To our knowledge, there is no study of such a generalization of Ramsey numbers in the literature.
Hence, we apply the probabilistic method to give a lower bound of $R^m(k;r)$.

\begin{lem}\label{gramseyl} Let $m,k,r$ be integers such that $2\leq m<r$ and $k>1$. Then $R^{m}(k;r)\geq \frac{\sqrt{m}}{\sqrt{e}r}k\left(\frac{r}{m}\right)^{k/2}$ provided that $\left(\frac{r}{m}\right)^{\frac{k-1}{2}}\left(\frac{e}{r}\right)^{\frac{1}{k}}\geq 2e$.
\end{lem}
\begin{pf}Colour the edges of the complete graph $K_n$ by ${[r]\choose m}$ randomly. That is, we
colour each edge by a set of $m$ colors with probability $\frac{1}{{r\choose m}}$.
Since the probability that a given copy of $K_k$ has all edges with a particular color $i$ is
$\left({r-1\choose m-1}/{r\choose m}\right)^{{k\choose 2}}=\left(m/r\right)^{{k\choose 2}}$,
 the
expected number of monochromatic copies of $K_k$ in color $i$ in this graph is
${n\choose k}\left(m/r\right)^{{k\choose 2}}$.
 Therefore, the
expected number of monochromatic copies of $K_k$ is
\[r{n\choose k}\left(\frac{m}{r}\right)^{{k\choose 2}}.\]
Substituting $n=\frac{\sqrt{m}}{e\sqrt{r}}k(\frac{r}{m})^{k/2}(\frac{e}{r})^{1/k}\geq 2k$, we have $r{n\choose k}\left(m/r\right)^{{k\choose 2}}<1$. Hence, we get \[R^{m}(k;r)\geq \frac{\sqrt{m}}{\sqrt{e}r}k\left(\frac{r}{m}\right)^{k/2}.\]
\qed \end{pf}

The construction below is a generalization of Construction~\ref{con4} by using extended Ramsey colorings.

\begin{cons}\label{con5} Let $t=2s+\delta$, $v=s+\alpha$, where $\delta=0$ or $1$, and $\alpha\geq 3$ is a constant. Let $c=t+2-v=s+\delta+2-\alpha$ and $r=v-c=2\alpha-\delta-2$. Let $R=[c+1,v]$, then $|R|=r$. Let $N=v(c-1)+l$ for some $l$, which will be determined later.  We construct an $N\times v$ array $C$ as follows.  For each $i\in [c]$ and each $j\in [v]$ and $j\neq i$, $ji$ starts a row. This accounts for $v(c-1)+r$ rows. Then let the symbol $c+i$ starts $k_i-1$ other rows, $k_i\geq r$, $i\in [r]$ and $k_1+\cdots+k_r=l$. Now we have in total $N$ rows. We can further assume that for each $i\in R$, and each different $j\in [v]$, $ij$ starts a row at least once.

Construct a complete graph $G=K_c$ with vertex set $[c]$. Suppose that  $c< R^{r-2}(k_1+1,\ldots,k_r+1)$, and we have a  Ramsey $(k_1+1,\ldots,k_r+1)^{r-2}$-coloring   of $K_c$. That is, there does not exist any monochromatic $K_{k_h+1}$ subgraph of color $h$ for $h \in [r]$.  Now for each pair $i,j\in [c]$, there are exactly two rows starting with $ij$ or $ji$. If the edge $i,j$ in $G$ is colored by a set $H$ of $r-2$ colors, then assign the third elements of these two rows by $c+h'$, where $h'\in [r]\setminus H$.  Finally, complete each row arbitrarily to a permutation over $[v]$.
\end{cons}

The proof of the following lemma is similar to that of Lemma~\ref{ramsey1} thus omitted here.

\begin{lem}\label{ramsey2} Let $t=2s+\delta$, $v=s+\alpha$, where $\delta=0$ or $1$, and $\alpha\geq 3$ is a constant. Let $c=s+\delta+2-\alpha$ and $r=2\alpha-\delta-2$. Let $l$ be a positive integer such that there exists $r$ integers $k_i\geq r$ and $k_1+\cdots+k_r=l$. Then there exists an  $(N,s+\alpha,2s+\delta)$-suitable core with $N=v(c-1)+l$ provided that $c< R^{r-2}(k_1+1, \ldots, k_r+1)$.
\end{lem}

Hence, by applying Lemma~\ref{ramsey2} with $k_i=l/r$, $i\in [r]$, we can prove Theorem~\ref{upbound2}.

%\begin{thm}\label{upbound2}  Let $t=2s+\delta$, $v=s+\alpha$, where $\delta=0$ or $1$, and $\alpha\geq 3$ be a constant.
% For each constant $\sigma\geq \frac{2d}{\ln (d)-\ln(d-2)}$ with $d=2\alpha-\delta-2$,  there exists an $(N,v,t)$-suitable core with $N=v(t+1-v)+\sigma\ln (s)$ for all large $s$.
%\end{thm}
\begin{pot5}
Let $c$ and $r$ be defined as in Lemma~\ref{ramsey2}.  Then there exists an  $(N,s+\alpha,2s+\delta)$-suitable core with $N=v(c-1)+l$ provided that $c< R^{r-2}(k_1+1, \ldots, k_r+1)$ for some integers $k_i\geq r$ and $k_1+\cdots+k_r=l$. Set  $k_i=l/r$, $i\in [r]$. By Lemma~\ref{gramseyl}, we have $R^{r-2}(l/r+1;r)\geq \left(\frac{r}{r-2}\right)^{l/(2r)}> s\geq c$ when $l\geq \frac{2r}{\ln r-\ln(r-2)}\ln s$.\qed
 \end{pot5}

To be complete,  we generalize the classic Erd\"{o}s-Szekeres argument  to get an upper  bound for the extended Ramsey number $R^m(k;r)$.

 \begin{lem}\label{gramsey2}$R^{m}(k_1,k_2,\ldots,k_r)\leq\frac{1}{m}(R_1+R_2+\cdots+R_r)$, where $R_i=R^{m}(k_1,k_2,\ldots,k_{i-1},k_i-1,k_{i+1},\ldots,k_r)$, $i\in[r]$.
\end{lem}
\begin{pf}Let $n=\frac{1}{m}(R_1+\cdots+R_r)$, and let $L$ be any $(r;m)$-coloring of $K_n$. Define $S_i=\{u\in V(K_n): i\in L(u,v)\}$, $i\in [r]$. Since each $|L(u,v)|=m$, we have \[|S_1|+\cdots+|S_r|=m(n-1)=R_1+\cdots+R_r-m.\]
Since $m<r$, we must have $|S_i|=R_i$ for some $i\in [r]$ by pigeonhole principle. Consider the complete subgraph $G$  of $K_n$ with vertex set $S_i$ and the edge coloring induced by $L$. Since $G$ has $R_i$ vertices, $G$ has a monochromatic $K_{k_j}$ in color $j$, for some $j\neq i$, or a monochromatic $K_{k_i-1}$ in color $i$. Including the vertex $v$, we conclude that any coloring of $K_n$ has a monochromatic $K_{k_j}$ in color $j$, for some $j\in [r]$.
\qed \end{pf}

From the symmetry of $R^{m}(k_1,k_2,\ldots,k_r)$, we can assume that $k_1\leq k_2\leq \cdots \leq k_r$. It is easy to see that $R^{m}(2,2,\ldots,2,k_{r-m+1},\ldots,k_r)=k_{r-m+1}$. Then by induction on $\sum_{i=1}^{r} k_i$, the following is an immediate consequence of Lemma~\ref{gramsey2}.

\begin{cor}\label{gramsey3}
$R^{m}(k_1,k_2,\ldots,k_r)\leq \left(\frac{1}{m}\right)^{\nu}\frac{(k_1+k_2+\cdots+k_r-r)!}{(k_1-1)!(k_2-1)!\cdots(k_r-1)!}$, where $\nu=\sum_{i=1}^{r-m} (k_i-2)$. In particular, $R^m(k;r)\leq \left(\frac{1}{m}\right)^{\nu} r^{rk}$.
\end{cor}

%\begin{lem}\label{lbound2}   There exists $s_0$ such that an  $(N,s+\alpha,2s+1)$-suitable core with $N=(s+\alpha)(s+2-\alpha)+\frac{\ln (s)}{2r\ln r}$ does not exist for all $s\geq s_0$, where $r=2\alpha-3$.
%\end{lem}
%\begin{pf}  Let  $l= \frac{\ln (s)}{2r\ln r}$. Set $d=3ls^{1/2}$ and $e=s^{1/2}$. The proof is similar to that of Lemma~\ref{lbound}, except the function $f$ maps $[c]$ to $r$-subsets of $R$, and all symbols in $B$ are mapped the same $r$-set $F=\{r_1,r_2,\ldots,r_r\}\subset R$. Then color the edges of the complete graph $G$ with vertex set $B$ by ${F\choose {r-2}}$  as follows: for each $b_i,b_j\in B$, there is at at least $r-2$ elements of $F$ that each do not precedes all $r-1$ others in neither the row starting $b_ib_j$ nor the row starting $b_jb_i$; choose $r-2$ such elements to form a set $H\subset F$ and color the edge between $b_i$ and $b_j$ by the set $H$. So the resulting graph $G$ is an $(r;r-2)$-coloring of $K_e$. Suppose in $C$, the symbols $r_i$, $i\in [r]$ start a row exactly $s+2-\alpha+k_i$ times, where $k_i\leq l-1$. Since $e\geq \left(\frac{1}{r-2}\right)^{2l-4}r^{rl}\geq  R^{r-2}(l;r)\geq R^{r-2}(k_1+1,\ldots,k_r+1)$, there does not exist a Ramsey $(k_1+1,\ldots,k_r+1)^{r-2}$-coloring of $K_e$, which enables us to find the contradiction. \qed
%\end{pf}

 \section{Concluding remark}
 \label{sec:con}
 In this paper, we give new existence and nonexistence proofs of $(N,v,t)$-suitable cores from packings of triples and Ramsey colorings. Our main results Theorems~\ref{upbound},~\ref{asym} and~\ref{upbound2}  suggest that: for any fixed integer $\alpha\geq 2$,  there may exist some constants $\sigma_{\alpha}$ and $\tau_{\alpha}$, such that if $\sigma_{\alpha} \ln s\leq l(s)\leq \tau_{\alpha} \ln s$, then SCN$(t,N)=s+{\alpha}$ for all sufficiently large $s$, where $t=2s+0$ or $1$, and $N=(s+{\alpha})(t+1-s-{\alpha})+l$.

For ${\alpha}=2$, Theorem~\ref{asym} shows that $\sigma_{\alpha}=0$ and $\tau_{\alpha}=\frac{1}{6\ln 3}$ satisfy the condition. What are the possible values of $\sigma_{\alpha}$ and $\tau_{\alpha}$ for ${\alpha}>2$? Theorem~\ref{upbound2} gives an example of $\sigma_{\alpha}$ for each $\alpha>2$. For $\tau_{\alpha}$, we can try similar arguments as in Lemma~\ref{lbound}. But using the rough upper bound of the extended Ramsey number $R^m(k;r)$ in Corollary~\ref{gramsey3}, we can only get a value of $\tau_{\alpha}$ less than $\sigma_{\alpha}$, which is useless.

Finally, we mention that, Balandraud et al. \cite{balandraud2018largest} determined the maximum size of {\it minimal} $2$-suitable sets of permutations, where no proper subsets are $2$-suitable. This problem arises in the determination of the Carath\'eodory numbers for certain
abstract convexity structures on the $(n - 1)$-dimensional real and integer vector spaces.  It would be interesting to consider such questions with the objective of determining or estimating the maximum cardinality of a minimal $t$-suitable set of permutations for $t\geq 3$.
%Constructions~\ref{con1} give very balanced suitable cores. Are they the largest minimal suitable cores, as studied in \cite{balandraud2018largest} for $t=2$?


\section*{Acknowledgments}
This research is supported by NSFC  under grants 11771419 and 11301503, and by ``the  Fundamental
Research Funds for the Central Universities''.

\bigskip\noindent{\bf References}

%\bibliographystyle{elsarticle-num}
%\bibliography{D:/work/1/JabRefdata/mine}

\end{document}